\documentclass[12pt]{amsart}
\usepackage[cp1251]{inputenc}
\usepackage[english]{babel}
\usepackage{amsmath}
\usepackage{amssymb}
\usepackage{amsfonts}

\begin{document}

\title[On asymptotical expansions for certain singular integrals]{On asymptotical expansions for certain singular integrals: 2-dimensional case}

\author[V.B. Vasilyev]{Vladimir B. Vasilyev}
\date{}

\maketitle

{\bf Abstract.} { One discusses a problem of asymptotical behavior for some operators in a general theory of pseudo differential equations on manifolds with borders. Using the distribution theory one obtains certain explicit representations for these operators. These limit distributions are constructed with a help of the Fourier transform, Dirac mass-function and its derivatives, and well-known distribution related to the Cauchy type integral.}



\section{Introduction}
\label{sec:1}

In the theory of pseudo differential equations the main difficulty is studying model operators in canonical domains according to a local principle. It declares that for a Fredholm property of a general pseudo differential operator on a compact manifold one needs an invertibility of its local representatives in each point of the manifold. The author wrote many times on a nature of these local representatives, these are distinct in dependence on a point of a manifold. Each singularity (half-space, cone, wedge, etc.) corresponds to a certain distribution, and a convolution with the distribution describes and something else describe a local representative of an initial pseudo differential operator in an appropriate point of manifold. All details can be found in author's studies \cite{V0,V1,V4}. But singularities can be of distinct dimensions and it is possible such singularities of a low dimension can be obtained from analogous singularities of full dimension. It means we need to find distributions for limit cases when some of parameters of singularities tend to zero. This approach was partially realized in author's papers \cite{V2,V3}, and the latest paper \cite{V4} is devoted to multi-dimensional constructions. The further author's idea is the following.  If we know the limit operator for a thin singularity then possible it is zero approximation for a such thin singularity. It is desirable to obtain an asymptotic expansion with a small parameter for the distribution corresponding to a such singularity. We will consider here a two-dimensional case and hope these studies will help us to transfer such constructions to a multi-dimensional situations \cite{V4}.

\section{An initial approximation}

To describe a solvability picture for a model elliptic pseudo differential equation with an operator $A$
\begin{equation}\label{1}
(Au)(x)=v(x),
\end{equation}
in multi-dimensional cone $C^a_+=\{x\in\mathbb R^2: x_2>a|x_1|,  a>0\}$
the author earlier considered a special singular integral operator \cite{V1}
\[
(K_au)(x)=\frac{a}{2\pi ^{2} }\lim_{\tau\to 0+}\int\limits_{\mathbb R^2}\frac{u(y)dy}{(x_1-y_1)^2-a^2(x_2-y_2+i\tau)^2}.
\]
This operator served a conical singularity in the general theory of boundary value problems for elliptic pseudo differential equations on manifolds with a non-smooth boundary. The operator \eqref{1} is a convolution operator, and the parameter $a$ is a size of an angle, $x_2>a|x_1|, a=\cot\alpha$.

We will consider two spaces of basic functions for distributions. If $D(\mathbb R^2)$ denotes a space of infinitely differentiable functions with a compact support then $D'(\mathbb R^2)$ is the corresponding space of distributions over the space $D(\mathbb R^2)$, analogously if $S(\mathbb R^2)$ is the Schwartz space of  infinitely differentiable rapidly decreasing at infinity functions then $S'(\mathbb R^2)$ is a corresponding space of distributions over $S(\mathbb R^2)$.

When $a\to +\infty$ one obtains \cite{V2} the following limit distribution
\begin{equation} \label{2}
 \mathop{\lim }\limits_{a\to \infty } \frac{a}{2\pi ^{2} } \frac{1}{\xi _{1}^{2} -a^{2} \xi _{2}^{2} } =\frac{i}{2\pi } {\mathcal P} \frac{1}{\xi _{1} } \otimes \delta \left(\xi _{2} \right),
  \end{equation}
where the notation for distribution $\mathcal P$  is taken from V.S. Vladimirov's books \cite{Vl1,Vl2}, and  $\otimes $
denotes the direct product of distributions. Here $\delta$ denotes one-dimensional Dirac mass-function which acts on $\varphi\in D(\mathbb R)$ by the following way
\[
(\delta, \varphi)=\varphi(0),
\]
and the distribution  $\mathcal P\frac{1}{x}$ is defined by the formula
\[
(\mathcal P\frac{1}{x}, \varphi)=v.p.\int\limits_{-\infty}^{+\infty}\frac{\varphi(x)dx}{x}\equiv\lim\limits_{\varepsilon\to 0+}\left(\int\limits_{-\infty}^{-\varepsilon}+\int\limits_{\varepsilon}^{+\infty}\right)\frac{\varphi(x)dx}{x}.
\]

Our main goal in this paper is obtaining an asymptotical expansion for the two-dimensional distribution
\[
K_a(\xi_1,\xi_2)\equiv \frac{a}{2\pi ^{2} } \frac{1}{\xi _{1}^{2} -a^{2} \xi _{2}^{2} }
\]
with respect to small $a^{-1}$. It is defined by the corresponding formula $\forall\varphi\in D(\mathbb R^2)$
\begin{equation}\label{3}
(K_a, \varphi)=\frac{a}{2\pi ^{2} }\int\limits_{\mathbb R^2}\frac{\varphi(\xi_1,\xi_2)d\xi}{\xi _{1}^{2} -a^{2} \xi _{2}^{2} }.
\end{equation}

\section{A decomposition formula for distributions}

We will use the standard Maclaurin formula
\[
\varphi(\xi_1,\xi_2)=\sum\limits_{k=0}^{\infty}\frac{\varphi^{(k)}_{\xi_2}(\xi_1,0)}{k!}\xi_2^k
\]
and make the change of a variable $a\xi_2=t, a^{-1}=b,$ then the formula \eqref{3} will become
\[
(K_a, \varphi)=\frac{1}{2\pi ^{2} }\int\limits_{\mathbb R^2}\frac{\varphi(\xi_1,bt)d\xi_1dt}{\xi _{1}^{2} -t^{2} }.
\]

Then we represent $\mathbb R^2=M\cup(\mathbb R^2\setminus M)$ where $M$ is a square with a line size $N$, so we have
\begin{equation}\label{4}
(K_a, \varphi)=\frac{1}{2\pi ^{2} }\left(\int\limits_{M}+\int\limits_{\mathbb R^2\setminus M}\right)\frac{\varphi(\xi_1,bt)d\xi_1dt}{\xi _{1}^{2} -t^{2} }
\end{equation}

\subsection{A rough decomposition}

Let us consider here $\varphi\in D(\mathbb R^2)$. Since a support of $\varphi$ is a compact set we have a one summand in the formula \eqref{4} therefore we might work with the formula
\[
\frac{1}{2\pi ^{2} }\int\limits_{\mathbb R^2}\frac{\varphi(\xi_1,bt)d\xi_1dt}{\xi _{1}^{2} -t^{2} }
\]
immediately.

More naturally for the last will be by the following way using Maclaurin series
\[
(K_a, \varphi)=\frac{1}{2\pi ^{2} }\sum\limits_{k=0}^{\infty}\frac{1}{k!}\int\limits_{\mathbb R^2}
\frac{\varphi^{(k)}_{\xi_2}(\xi_1,0)b^kt^kd\xi_1dt}{\xi_1^2-t^2}
\]

If $t$ varies in a line segment then $bt\sim b, b\to 0$, and we can use
the following formal representations \cite{V3}
$$
K_a(\xi_1,\xi_2)=\frac{i}{2\pi }\sum\limits_{n=0}^{+\infty}\frac{(-1)^n}{n!a^n}{\mathcal P}\frac{1}{\xi_1}\otimes\delta^{(n)}(\xi_2).
$$

\subsection{A sharp decomposition}

Here we consider $\varphi\in S(\mathbb R^2)$.

A formal using the Maclaurin formula for the first integral in \eqref{4} will lead to the following result
\begin{equation}\label{5}
(K_a, \varphi)=\frac{1}{2\pi ^{2} }\sum\limits_{k=0}^{\infty}\frac{b^k}{k!}\int\limits_{-N}^{+N}\varphi^{(k)}_{\xi_2}(\xi_1,0)
\left(\int\limits_{-N}^{+N}\frac{t^kdt}{\xi_1^2-t^2}\right)d\xi_1,
\end{equation}
and we need to give a certain sense for the expression in brackets.

Let us denote
\[
T_{k,N}(\xi_1)\equiv\int\limits_{-N}^{+N}\frac{t^kdt}{\xi_1^2-t^2}
\]
and reproduce some calculations.

First $T_{k,N}(\xi_1)\equiv 0, \forall k=2n-1, n\in\mathbb N$. So the non-trivial case is $k=2n, n\in\mathbb N$. Let us remind $T_{0,\infty}(\xi_1)=\pi i2^{-1}\xi_1^{-1}$ \cite{V2,V3}. For other cases we can calculate this integral. we have the following

$k=2$,
$$
T_{2,N}(\xi_1)=-2N-2^{-1}\xi_1^{-1}\ln\frac{N-\xi_1}{N+\xi_1}+\pi i2^{-1}\xi_1^{-1};
$$

$k=4$,
$$
T_{4,N}(\xi_1)=-2/3N^3-2\xi_1^2N-2^{-1}\xi_1^3\ln\frac{N-\xi_1}{N+\xi_1}+\pi i2^{-1}\xi_1^{3};
$$

$k=6$,
$$
T_{6,N}(\xi_1)=-2/5N^5-2/3\xi_1^2N^3-2\xi_1^5N-2^{-1}\xi_1^5\ln\frac{N-\xi_1}{N+\xi_1}+\pi i2^{-1}\xi_1^{5},
$$
and so on. One can easily write all expressions for arbitrary $T_{2n,N}(\xi_1)$.

 In general one can write
 $$
 T_{2n,N}(\xi_1)=P_{2n-1}(N,\xi_1)-2^{-1}\xi_1^{2n-1}\ln\frac{N-\xi_1}{N+\xi_1}+\pi i2^{-1}\xi_1^{2n-1}
 $$
where $P_{2n-1}(N,\xi_1)$ is a certain polynomial of order $2n-1$ on variables $N,\xi_1$.

Therefore instead of the formula \eqref{5} we can write
\begin{equation}\label{6}
(K_a, \varphi)=\frac{i}{2\pi}(\mathcal P\frac{1}{\xi_1}\otimes\delta(\xi_2),\varphi)+
\end{equation}
\[
\frac{1}{2\pi ^{2} }\sum\limits_{n=1}^{\infty}\frac{b^{2n}}{{(2n)}!}\int\limits_{-N}^{+N}\varphi^{(2n)}_{\xi_2}(\xi_1,0)
\left(P_{2n-1}(N,\xi_1)-2^{-1}\xi_1^{2n-1}\ln\frac{N-\xi_1}{N+\xi_1}+\pi i2^{-1}\xi_1^{2n-1}\right)d\xi_1.
\]

Let us describe the polynomial $P_{2n-1}(N,\xi_1)$ more precisely. Obviously
\[
P_{2n-1}(N,\xi_1)=c_{2n-1}N^{2n-1}+c_{2n-3}N^{2n-3}\xi_1^2+\cdots+c_1N\xi_1^{2n-1}.
\]

Further we rewrite the equality \eqref{6} in the following form
\[
(K_a, \varphi)=\frac{i}{2\pi}(\mathcal P\frac{1}{\xi_1}\otimes\delta(\xi_2),\varphi)+
\]
\[
\frac{1}{2\pi ^{2} }\sum\limits_{n=1}^{\infty}\frac{b^{2n}}{{(2n)}!}\sum\limits_{k=1}^nc_{2k-1}N^{2k-1}\int\limits_{-N}^{+N}\varphi^{(2n)}_{\xi_2}(\xi_1,0)\xi_1^{2k-1}d\xi_1-
\]
\[
\frac{1}{4\pi ^{2} }\sum\limits_{n=1}^{\infty}\frac{b^{2n}}{{(2n)}!}\int\limits_{-N}^{+N}\varphi^{(2n)}_{\xi_2}(\xi_1,0)\xi_1^{2n-1}\ln\frac{N-\xi_1}{N+\xi_1}d\xi_1+\frac{i}{4\pi  }\sum\limits_{n=1}^{\infty}\frac{b^{2n}}{{(2n)}!}\int\limits_{-N}^{+N}\varphi^{(2n)}_{\xi_2}(\xi_1,0)\xi_1^{2n-1}d\xi_1
\]

We will start from two last summands. The second summand does not play any role because
\[
\lim_{N\to+\infty}\ln\frac{N-\xi_1}{N+\xi_1}=0.
\]

The third summand we will represent according to lemma 1 (see below) taking into account that we can pass to the limit under $N\to+\infty$
\[
\frac{i}{4\pi  }\sum\limits_{n=1}^{\infty}\frac{b^{2n}}{{(2n)}!}(\widetilde{\delta^{(2n-1)}}(\xi_1)\otimes\delta^{(2n)}(\xi_2), \varphi)
\]

For the first summand we consider separately the case $Nb\sim 1 (N\to\infty, b\to 0)$. in other words we consider a special limit to justify the decomposition. Then
\[
\frac{1}{2\pi ^{2} }\sum\limits_{n=1}^{\infty}\frac{b^{2n}}{{(2n)}!}\sum\limits_{k=1}^nc_{2k-1}N^{2k-1}\int\limits_{-N}^{+N}\varphi^{(2n)}_{\xi_2}(\xi_1,0)\xi_1^{2k-1}d\xi_1\sim
\]
\[
\frac{1}{2\pi ^{2} }\sum\limits_{n=1}^{\infty}\frac{1}{{(2n)}!}\sum\limits_{k=1}^nc_{2k-1}b^{2n-2k+1}\int\limits_{-\infty}^{+\infty}\varphi^{(2n)}_{\xi_2}(\xi_1,0)\xi_1^{2k-1}d\xi_1.
\]

Therefore
\[
\frac{1}{2\pi ^{2} }\sum\limits_{n=1}^{\infty}\frac{1}{{(2n)}!}\sum\limits_{k=1}^nc_{2k-1}b^{2n-2k+1}\int\limits_{-\infty}^{+\infty}\varphi^{(2n)}_{\xi_2}(\xi_1,0)\xi_1^{2k-1}d\xi_1=
\]
\[
\frac{1}{2\pi^2 }\sum\limits_{n=1}^{\infty}\frac{1}{{(2n)}!}\sum\limits_{k=1}^nc_{2k-1}b^{2n-2k+1}(\widetilde{\delta^{(2k-1)}}(\xi_1)\otimes\delta^{(2n)}(\xi_2), \varphi)
\]

One can note if desirable
\[
c_{2k-1}=-2(1+\frac{1}{3}+\cdots+\frac{1}{2k-1}).
\]

\section{A local behavior of a boundary operator}

{\bf Lemma 1.}{\it 
If a distribution $a$ acts on the function $\varphi\in S(\mathbb R)$ by the following way
\[
(a,\varphi)=\int\limits_{-\infty}^{+\infty}\xi^k\varphi(\xi)d\xi,
\]
then this distribution $a$ is the following
\[
a(\xi)=\widetilde{\delta^{(k)}}(\xi),
\]
where the sign $\sim$ means here the inverse Fourier transform $F^{-1}$.
}

{\bf Proof.} Indeed we have $F\delta={\bf 1}$, where ${\bf 1}$ is an identity in a distribution sense so that $F^{-1}{\bf 1}=\delta$. Since
\[
(F(\varphi^{(k)}))(\xi)=(-1)^k\xi^k\tilde\varphi(\xi)
\]
then denoting $\psi=F^{-1}\varphi$ we write
\[
(a,\varphi)=(a,F\psi)=\int\limits_{-\infty}^{+\infty}\xi^k\widetilde\psi(\xi)d\xi=({\bf 1}, \xi^k\widetilde\psi(\xi))=({\bf 1},FF^{-1}(\xi_k\widetilde\psi(\xi)))=
\]
\[
(F{\bf 1},F^{-1}(\xi_k\widetilde\psi(\xi)))=(F{\bf 1},(-1)^k\psi^{(k)}(x))=(\delta,(-1)^k\psi^{(k)}(x))=
\]
\[
(\delta^{(k)},\psi)=(\delta^{(k)},F^{-1}\varphi)=(F^{-1}\delta^{(k)},\varphi),
\]
so we have the required identity.
$\triangle$

{\bf Theorem 1.}{\it 
The following formula
\[
K_a(\xi_1,\xi_2)=\frac{i}{2\pi}\mathcal P\frac{1}{\xi_1}\otimes\delta(\xi_2)+\sum\limits_{m,n}c_{m,n}(a)\widetilde{\delta^{(m)}}(\xi_1)\otimes\delta^{(n)}(\xi_2),
\]
where $c_{m,n}(a)\to 0, a\to +\infty,$ holds in a distribution sense.
}

{\bf Proof.} Returning to the formula \eqref{5} and using calculations $T_{k,N}(\xi_1)$ and lemma 1 we obtain the required assertion.
$\triangle$

{\bf Remark 1.} {\it
One can easily reconstruct coefficients $c_{m,n}(a)$ starting from above calculations.
}

\section{Towards a pseudo differential equation}

Let us return to the equation \eqref{1}. We will remind some our preliminary results \cite{V0,V1}.
The symbol $\stackrel{*} {C^a_+}$ denotes a conjugate cone for $C^a_+$:
\[
\stackrel{*} {C^a_+}=\{x\in{\mathbb R}^2: x=(x_1,x_2), ax_2>|x_1|\},
\]
$C^a_-\equiv -C^a_+,~T(C^a_+)$ denotes a radial tube domain over the cone $C^a_+$ \cite{Vl2}, i.e. domain in a complex space ${\mathbb C}^2$ of type ${\mathbb R}^2+iC^a_+$.

We consider symbols $A(\xi)$ satisfying the condition
$$
c_1\leq|A(\xi)(1+|\xi|)^{-\alpha}|\leq c_2,
$$
which are elliptic, and the number $\alpha\in{\mathbb R}$ is called an order of the operator $A$.

To describe the solvability picture for the equation \eqref{1} we use the following

{\bf Definition.} {\it
Wave factorization with respect to the cone $C^a_+$ for the symbol $A(\xi)$ is called its representation in the form
$$
A(\xi)=A_{\neq}(\xi)A_=(\xi),
$$
where the factors $A_{\neq}(\xi),A_=(\xi)$ must satisfy the following conditions:

1) $A_{\neq}(\xi),A_=(\xi)$ are defined for all admissible values $\xi\in{\mathbb R}^2$, without may be, the points $\{\xi\in{\mathbb R}^2:|\xi_1|^2=a^2\xi^2_2\}$;

2) $A_{\neq}(\xi),A_=(\xi)$ admit an analytical continuation into radial tube domains $T(\stackrel{*} {C^a_+}),T(\stackrel{*} {C^a_-})$ respectively with estimates
$$
|A_{\neq}^{\pm 1}(\xi+i\tau)|\leq c_1(1+|\xi|+|\tau|)^{\pm\kappa},
$$
$$
|A_{=}^{\pm 1}(\xi-i\tau)|\leq c_2(1+|\xi|+|\tau|)^{\pm(\alpha-\kappa)},~\forall\tau\in\stackrel{*} {C^a_+}.
$$

The number $\kappa\in{\mathbb R}$ is called index of wave factorization.
}

For $|\kappa-s|<1/2$ one has the existence and uniqueness theorem \cite{V0}
\[
\widetilde u(\xi)=A_{\neq}^{-1}(\xi)(K_a\widetilde {lv})(\xi),
\]
where $lv$ is an arbitrary continuation of $v$ on the whole $H^s(\mathbb R^2)$.

Below we denote $lv\equiv V$.

{\bf Theorem 2.} {\it
If the symbol $A(\xi)$ admits a wave factorization with respect to the cone $C^a_+$ and $|\kappa-s|<1/2$ then the equation $\eqref{1}$ has a unique solution in the space $H^s(C^a_+)$, and for a large $a$ it can be represented in the form
\[
\widetilde u(\xi)=\frac{i}{2\pi}A_{\neq}^{-1}(\xi)v.p.\int\limits_{-\infty}^{+\infty}\frac{(A_{=}^{-1}\widetilde V)(\eta_1,\xi_2)d\eta_1}{\xi_1-\eta_1}+
\]
\begin{equation}\label{7}
A_{\neq}^{-1}(\xi)\sum\limits_{m,n}c_{m,n}(a)\int\limits_{-\infty}^{+\infty}(\xi_1-\eta_1)^m(A_{=}^{-1}\widetilde V)^{(n)}_{\xi_2}(\eta_1,\xi_2)d\eta_1
\end{equation}
assuming $\widetilde{V}\in S(\mathbb R^2), A_{=}^{-1}\widetilde V$ means the function $A_{=}^{-1}(\xi)\widetilde V(\xi)$.
}

{\bf Proof.} We need to apply the theorem 1 and to recall correlations between distributions and pseudo differential operators. It proves the theorem.
$\triangle$

{\bf Remark 2.} {\it
A reader can easily write an analogue of the theorem 5.2 corresponding to a rough decomposition.
}

\section*{Conclusion}

It was shown the solution of the equation \eqref{1} for enough smooth right hand side $v$ can be represented in the form \eqref{7}.
It shows that in this series the first summand belongs to the space $H^s(C^a_+)$ only. Secondary summands can be useful for certain special situations related to some additional properties of the right hand side $v$.

\vspace{5cm}

Vladimir B. Vasilyev\\
Chair of Pure Mathematics\\
Lipetsk State Technical University\\
Lipetsk 398600, Russia\\
e-mail: vbv57@inbox.ru

\end{document}